\documentclass[oneside, reqno]{amsart}

\usepackage{xypic}
\input xy
\xyoption{all}
\usepackage{epsfig}
\usepackage{amsthm}
\usepackage{amssymb}
\usepackage{amsmath}
\usepackage{amscd}
\usepackage{color}
\usepackage[T1]{fontenc}
\usepackage{stackengine}
\usepackage{accents}
\stackMath

\newcommand{\bg}{\begin{equation}}
\newcommand{\ed}{\end{equation}}
\newcommand{\bga}{\begin{eqnarray}}
\newcommand{\eda}{\end{eqnarray}}
\newcommand{\pf}{\textbf{Proof:\ }}

\def\cbdu{\par{\raggedleft$\Box$\par}}

\newtheorem {Theorem}  {Theorem}

\numberwithin{Theorem}{section}

\newtheorem {Lemma}[Theorem]  {Lemma}
\newtheorem {Proposition}[Theorem]{Proposition}
\theoremstyle{definition}
\newtheorem{Definition}[Theorem]{Definition}
\theoremstyle{remark}

\expandafter\chardef\csname pre amssym.def
at\endcsname=\the\catcode`\@ \catcode`\@=11
\def\undefine#1{\let#1\undefined}
\def\newsymbol#1#2#3#4#5{\let\next@\relax
 \ifnum#2=\@ne\let\next@\msafam@\else
 \ifnum#2=\tw@\let\next@\msbfam@\fi\fi
 \mathchardef#1="#3\next@#4#5}
\def\mathhexbox@#1#2#3{\relax
 \ifmmode\mathpalette{}{\m@th\mathchar"#1#2#3}%
 \else\leavevmode\hbox{$\m@th\mathchar"#1#2#3$}\fi}
\def\hexnumber@#1{\ifcase#1 0\or 1\or 2\or 3\or 4\or 5\or 6\or 7\or 8\or
 9\or A\or B\or C\or D\or E\or F\fi}

\font\teneufm=eufm10 \font\seveneufm=eufm7 \font\fiveeufm=eufm5
\newfam\eufmfam
\textfont\eufmfam=\teneufm \scriptfont\eufmfam=\seveneufm
\scriptscriptfont\eufmfam=\fiveeufm

\catcode`\@=\csname pre amssym.def at\endcsname

\newcounter{remark}
\def  \T   {{\mathbb T}}

\def  \12  {{\frac{1}{2}}}

\def\build#1_#2^#3{\mathrel{\mathop{\kern 0pt#1}\limits_{#2}^{#3}}}

\begin{document}

\title[Anomalous dissipation for electron-MHD]{Anomalous dissipation of energy and magnetic helicity for the electron-MHD system}

\author [Mimi Dai]{Mimi Dai}

\address{Department of Mathematics, Stat. and Comp.Sci., University of Illinois Chicago, Chicago, IL 60607, USA}
\email{mdai@uic.edu} 

\author [Han Liu]{Han Liu}

\address{Department of Mathematics, Stat. and Comp.Sci., University of Illinois Chicago, Chicago, IL 60607, USA}
\email{hliu94@uic.edu}

\thanks{The authors were partially supported by NSF grant
DMS--1815069.}

\begin{abstract}

Via convex integration, we construct weak solutions to the non-resistive electron-MHD system which do not conserve energy and magnetic helicity.

\bigskip


\hspace{0.02cm}CLASSIFICATION CODE: 35Q35, 35D30, 76W05.
\end{abstract}

\maketitle

\section{Introduction}

The incompressible Hall-magneto-hydrodynamics (Hall-MHD) system
\begin{equation}\label{HMHD}
\begin{cases}
u_t+u\cdot\nabla u-B\cdot\nabla B+\nabla p= \nu\Delta u,\\
B_t+u\cdot\nabla B-B\cdot\nabla u+d_i\nabla\times((\nabla\times B)\times B)= \mu\Delta B,\\
\nabla \cdot u= 0,\\
u(0,x)=u_0(x),\  B(0, x)=B_0(x), \ t \in \mathbb{R}^+, \ x \in \mathbb{R}^3 \text{ or } \mathbb{T}^3
\end{cases}
\end{equation}
has recently drawn mathematicians' attentions, owing to its significance in plasma physics as well as its mathematical complexity. In system \eqref{HMHD}, the vector fields $u$ and $B$ are the fluid velocity and magnetic field, respectively, while $p$ denotes the scalar pressure. The coefficients $\nu,\mu$ and $d_i$ represent the kinematic viscosity, magnetic resistivity and ion inertial length, respectively. Given $B_0$ with $\nabla\cdot B_0=0,$ we can easily deduce that $\nabla\cdot B=0$ for all $t\geq 0.$ Hence, we assume that $\nabla \cdot B =0$ throughout the paper.

For the ideal MHD system, i.e., $\nu=\mu=d_i=0$ in \eqref{HMHD}, Kelvin's circulation theorem for any material surface $\mathcal{S}$ moving with the MHD fluid yields 
$$\frac{d}{dt}\int_{\mathcal S}B\cdot \mathrm{d}S=0,$$ that is to say, the magnetic flux through any material surface advected by the fluid is conserved. This is known as Alfv\'en's theorem, indicating that magnetic field lines are ``frozen'' into the fluid and the magnetic topology is preserved. Yet, for $d_i>0$, the Hall effect breaks the ``frozen-in'' property. Moreover, at length scale $\ell \ll d_i$ the ions appear to be too heavy to move with the fluid  and decouple from the magnetic field, leaving the magnetic field frozen into the electronic fluid only. In this situation, the background ionic flow velocity $u$ vanishes, and the Hall-MHD system \eqref{HMHD} reduces to the following electron-MHD (EMHD) system
\begin{equation}\label{emhd}
\begin{cases}
B_t+d_i\nabla\times((\nabla\times B)\times B)=\mu\Delta B,\\
\nabla \cdot B =0,\\
B(0, x)=B_0(x), \ t \in \mathbb{R}^+, \ x \in \mathbb{R}^3 \text{ or } \mathbb{T}^3.
\end{cases}
\end{equation}

In the past decade, various mathematical results on the Hall-MHD and EMHD systems have been obtained. Global existence of weak solutions was shown in \cite{ADFL, CDL, DS}, while well-posedness in various spaces was established in \cite{CDL, CWW, D1, D2, DL2}. Striking ill-posedness results can be found in \cite{CW, JO}, whereas non-uniqueness of weak solutions was proven in \cite{D3} via a convex integration scheme. In \cite{CS, DL}, asymptotic behavior of solutions was studied. For a variety of regularity and blow-up criteria, we refer readers to \cite{CL, CW, D4}.

We assume that there is a magnetic vector potential by $A$ satisfying $B=\nabla \times A,$ which can be chosen to be divergence-free under the assumption of Coulomb gauge. System (\ref{emhd}) implies that $A$ satisfies the following system of equations
\begin{equation}\label{amhd}
\begin{cases}
A_t+d_i(\nabla\times B)\times B=\mu\Delta A, \\
\nabla \cdot B =0, \\
B(0, x)=B_0(x), \ t \in \mathbb{R}^+, \ x \in \mathbb{R}^3 \text{ or } \mathbb{T}^3.
\end{cases}
\end{equation} 
Given the magnetic field $B,$ we can recover the magnetic potential $A$ through the Biot-Savart law, i.e., $$A=\nabla \times (-\Delta)^{-1} B.$$ 

In this paper, we will consider system \eqref{amhd} in the non-resistive setting $\mu =0$. We are interested in the magnetic energy $\mathcal E(t) = \|B(t)\|_{L^2}$ and magnetic helicity $\mathcal H(t) = \langle A, B\rangle_{L^2}(t),$ which are conserved over the time for sufficiently regular solutions. Indeed, the following identities hold --
\begin{equation}\notag
\frac12\frac{\mathrm d}{\mathrm dt}\mathcal E(t)+\int_{\mathbb T^3} \nabla\times \big((\nabla\times B)\times B \big)\cdot B\, \mathrm{d}x=0,\\
\end{equation}
\begin{equation}\notag
\frac{\mathrm d}{\mathrm dt}\mathcal H(t)+\int_{\mathbb T^3} \big((\nabla\times B)\times B\big) \cdot B\, \mathrm dx+\int_{\mathbb T^3} \nabla\times((\nabla\times B)\times B)\cdot A\, \mathrm{d}x=0.
\end{equation}
It's easy to see that $\mathcal E'(t) =0$ and $\mathcal H'(t) =0,$ as  integration by parts along with the properties of cross product yield
\begin{gather*}
\int_{\mathbb T^3} \nabla\times((\nabla\times B)\times B)\cdot B\mathrm dx=\int_{\mathbb T^3} (\nabla\times B)\times B\cdot(\nabla\times B)\,\mathrm dx=0,\\
\int_{\mathbb T^3} \big((\nabla\times B)\times B\big) \cdot B\, \mathrm dx=0,\\
\int_{\mathbb T^3} \nabla\times((\nabla\times B)\times B)\cdot A\mathrm dx
= \int_{\mathbb T^3} (\nabla\times B)\times B\cdot (\nabla\times A)\, \mathrm dx = 0.
\end{gather*}

However, for solutions that are not regular enough, the operation of integration by parts is not justified and the aforementioned conservation laws may no longer be true as the otherwise conserved quantities seem to dissipate due to the lack of smoothness of the solutions. This mechanism, known as anomalous dissipation, lies at the core of Onsager's conjecture in the context of hydrodynamics, which predicts that H\"older regularity index $1/3$ is a threshold determining whether or not a solution to the Euler equations conserves energy. As of now, both the positive and negative sides of Onsager's conjecture have been mathematically justified by a series of studies on the 3D Navier- Stokes and Euler equations. The review articles \cite{BV, S} are great resources on the development of this topic and related ones.

A heuristic scaling analysis of the Hall term leads to the expectation that solutions belonging to the Onsager-critical spaces $ L^3(0,T; B^{2/3}_{3, c_0})$ and $L^3(0,T; B^{1/3}_{3, c_0})$ conserve the energy and the magnetic helicity, respectively. (The space $B^s_{3, c_0}$ is the union of all Besov spaces $B^s_{3, q}$ with finite summability index $q.$) Indeed, this was proven for system \eqref{HMHD} with $\nu=\mu=0$ in \cite{DS}. As the Hall term is the most singular term in system \eqref{HMHD} which breaks the symmetries of the MHD system, the properties of the Hall-MHD system can be inferred from those of the EMHD system. Hence, this paper deals with the EMHD system only and the study of the full Hall-MHD system shall be forthcoming . 

In this paper, using convex integration, we aim to construct weak solutions to the non-resistive EMHD system on $\mathbb T^3$ such that while the solutions have finite energies, the corresponding $\mathcal E(t)$ and $\mathcal H(t)$ are non-constant functions. 
Our result of non-conservation may serve as an interpretation of the magnetic reconnection phenomenon, often observed in space plasmas. As mentioned before, the Hall-MHD system is a pivotal model for magnetic reconnection, a process in which Alfv\'en's ``frozen-in'' theorem is violated and changes in the magnetic topology occur, accompanied by energy transfers. As discovered in \cite{M}, the magnetic helicity turns out to measure the degree of self-linkage or knottedness of the magnetic field lines. Thus, non-conservation of the magnetic helicity implies topological changes of the magnetic field, which is a key evidence of magnetic reconnection. Our main result states as follows.
\begin{Theorem}\label{thm2}
There exist weak solutions $B\in C([0,1]; H^\alpha(\mathbb T^3))$ to the ideal EMHD system, i.e., system (\ref{emhd}) with $\mu=0,$ for some small constant $\alpha>0$, such that the magnetic energy $\mathcal E(t)$ and magnetic helicity $\mathcal H(t)$ are non-constant functions on the time interval $[0,1]$.
\end{Theorem}

We would like to mention that the existence of anomalous dissipation of magnetic helicity for the 3D ideal MHD, i.e., system (\ref{HMHD}) with $\nu=\mu=d_i=0,$ was shown for weak solutions with finite energy in \cite{BBV}, which also demonstrates that there exist finite energy weak solutions to the ideal MHD that cannot be attained in the vanishing viscosity and resistivity limit. Moreover, in a very recent paper \cite{FLS}, bounded and compactly supported solutions to the ideal MHD system that do not conserve the total energy and cross helicity have been constructed. We refer readers to the expository articles \cite{BV, DS3} for the development and applications of convex integration.

We shall introduce the so-called convex integration method, which is essentially an iterative scheme outlined as follows. Let $q \in \mathbb N.$ We approximate system \eqref{amhd} with $\mu=0$ successively with the following system
\begin{equation}\label{Aq}
\begin{cases}
\partial_t A_q+\nabla\cdot (B_q\otimes B_q) -\nabla p_q =\nabla\cdot R_q,\\
\nabla\times A_q=B_q,\\
\nabla\cdot A_q=\nabla\cdot B_q=0,
\end{cases}
\end{equation}
with $p_q = \frac{|B_q|^2}{2}$ and $R_q$ being a certain symmetric traceless tensor determined by its predecessor $R_{q-1}.$ At each step of the iteration, we will construct a solution $(A_q, B_q)$ to system (\ref{Aq}) by adding a perturbation $(v_q, w_q),$ based on the Beltrami waves, to the preceding duo $(A_{q-1}, B_{q-1}).$ As the iterative algorithm runs, $\|R_q\|_{L^p} \to 0$ for some $p$ close to $1.$  Hence, $\left\{ (A_q, B_q) \right\}_{q \geq 1}$ converges in $\mathcal{D}'$ to some limit $(A, B),$ which is a weak solution to system \eqref{amhd}, i.e., 
\begin{equation}\notag
\begin{cases}
A_t +\nabla\cdot (B\otimes B)-\nabla \frac{|B|^2}{2}=0,\\
B=\nabla\times A, \\
\nabla\cdot A=\nabla\cdot B=0.
\end{cases}
\end{equation}
We then show that for such $(A, B),$ the corresponding $\mathcal E(t)$ and $\mathcal H(t)$ are non-constant on the time interval $[0,1]$.

Our proof of Theorem \ref{thm2} relies on the following proposition on the iterative procedure, which we shall prove in Section \ref{sec-prop}.
\begin{Proposition}\label{prop}
For a triple $(A_{q}, B_{q}, R_{q})$ that solves the approximating system (\ref{Aq}) such that
\begin{equation}\label{est-induction}
\|B_q\|_{L^2}\leq 1-\delta_q^{\frac12}, \ \ \|B_q\|_{C^1_{x,t}}\leq \lambda_q^2, \ \ \|R_q\|_{L^1}\leq c \delta_{q+1},
\end{equation}
there exists a triple $(A_{q+1}, B_{q+1}, R_{q+1})$ satisfying system (\ref{Aq}) and estimates (\ref{est-induction}) at the $(q+1)$-th level. Moreover, for $\alpha<\frac{\beta}{1+\beta}$ and a small constant $\varepsilon_b>0$, we have
\begin{equation}\label{est-increm}
\|B_{q+1}-B_q\|_{L^2}\leq \delta_{q+1}^{\frac12}, \ \ \|B_{q+1}-B_q\|_{H^\alpha}\leq \delta_{q+1}^{\varepsilon_b}.
\end{equation}
\end{Proposition}

In the above proposition, a frequency parameter $\lambda_q$ and an amplitude parameter $\delta_q$ are introduced in order to measure the stress tensor $R_q.$ More specifically, we define 
\begin{equation}\label{pp1}
\lambda_q:=a^{b^q}, \ \ \delta_q:=\lambda_q^{-2\beta},
\end{equation}
for some positive regularity parameter $\beta \ll 1$ and some large $a, b \in \mathbb N.$ The choices of the parameters $a, b$ and $\beta$ shall be made more precise in the following sections.

The rest of the paper is organized as follows. In Section \ref{sec-pre}, we fix the notations and present several preliminary lemmas. In Section \ref{pt}, we deduce Theorem \ref{thm2} from Proposition \ref{prop}. Sections \ref{sec-convex} and \ref{sec-prop} shall be dedicated to the construction of $(v_q, w_q)$ and the proof of Proposition \ref{prop}.

\bigskip

\section{Preliminaries}
\label{sec-pre}

\subsection{Notation}

We specify a few notations for the sake of brevity. We denote by $A \lesssim B$ an estimate of the form $A\leq C B$ with a certain constant $C$ and by $A\sim B$ an estimate of the form $C_1 B\leq A\leq C_2 B$ with certain constants $C_1$ and $C_2$. 

The space of symmetric $3 \times 3$ matrices is written as $\mathcal{S}^{3 \times 3},$ while $\mathcal{S}_0^{3 \times 3}$ denotes the space of symmetric traceless $3 \times 3$ matrices. We denote the trace-free part of the tensor product by $\mathring \otimes.$

We define the norm $\|f\|_{C^1_{x,t}}:=\|f\|_{L^\infty} + \|\nabla f\|_{L^{\infty}} + \|\partial_t f \|_{L^\infty}.$ 

The operator $\mathbb P_{\leq N}$ (or $\mathbb P_{\geq N}$) can be understood as the projection onto frequencies no higher (or no lower, respectively) than $N.$ In turn, $\mathbb P_{\neq 0}$ projects onto non-zero frequencies.

\subsection{Definition of weak solutions}

We define the weak solutions to the electron-MHD system.
\begin{Definition}\label{weak-sol}
A divergence-free vector field $B$ is said to be a weak solution to system \eqref{emhd} on $[0, T]$ if it satisfies the following integral equation
\begin{equation}\notag
\int_0^T\int_{\mathbb T^3}B\cdot \varphi_t+(B\otimes B):\nabla\nabla\times \varphi\,\mathrm dx\,\mathrm dt=\mu \int_0^T\int_{\mathbb T^3}\nabla B:\nabla \varphi\,\mathrm dx\,\mathrm dt
\end{equation}
for any $\varphi\in C^\infty_c([0,T]\times \mathbb T^3).$

Equivalently, a pair of divergence-free vector fields $(A, B)$ is said to be a weak solution to system \eqref{amhd}, if for any $\varphi\in C^\infty_c([0,T]\times \mathbb T^3),$ the following integral equations hold --
\begin{gather*}
\int_0^T\int_{\mathbb T^3}A\cdot \varphi_t+(B\otimes B):\nabla\varphi\,\mathrm dx\mathrm dt=\int_0^T\int_{\mathbb T^3}\nabla A:\nabla \varphi\,\mathrm dx\mathrm dt,\\
\int_0^T\int_{\mathbb T^3} A \times \nabla \varphi \,\mathrm{d}x \mathrm{d}t =\mu \int^T_0 \int_{\mathbb T^3} B \cdot \varphi \,  \mathrm{d}x \mathrm{d}t.
\end{gather*}
\end{Definition}

For $\mu >0,$ the existence of weak solutions can be shown via the standard Galerkin approximation procedure (cf. \cite{CDL}), while in the case of ideal electron-MHD system, the existence of weak solutions in general remains an unresolved issue.

\subsection{Geometrical preliminaries}\label{G}

We introduce the Beltrami waves, which are stationary solutions to the Euler equations (thus also to the electron-MHD system). Let $\Lambda$ be a finite subset of $\mathbb S^2\cap\mathbb Q^3$ with $\Lambda=-\Lambda$. If $k \in \Lambda,$ then $\{ k, -k \} \subset \Lambda.$ Given $\{ k, -k \} \subset \Lambda,$ we fix $k_1 \in \mathbb S^2\cap\mathbb Q^3$ such that $k_1 \perp k.$  We supplement $\{ k, k_1 \}$ with $k_2 \in \mathbb S^2\cap\mathbb Q^3$ such that $k_2 \perp \text{span}\{k, k_1\}.$ 

Associating with each $k \in \Lambda$ a constant $a_k\in \mathbb C$ satisfying $\bar a_k=a_{-k},$ we define the real-valued Beltrami wave as
$$W(x):=\sum_{k\in \Lambda}a_k W_k(x),$$
where $$W_k(x)=\frac1{\sqrt 2}(k_1+ik_2)e^{i\lambda k\cdot x}$$ with $\lambda$ a certain large integer for which $\lambda \Lambda\subset \mathbb Z^3.$ 

We note that the space of Beltrami waves contains linear spaces of fairly large dimensions. The Beltrami wave possesses the following properties -- 
\begin{gather*}
\nabla\cdot W=0,  \ \  \nabla\times W=\lambda W,\\
\nabla\cdot (W\otimes W)=\nabla\frac{|W|^2}{2},
\end{gather*}
\[\displaystyle\stackinset{c}{}{c}{}{-\mkern4mu}{\displaystyle\int_{\mathbb T^3}} W\otimes W\, \mathrm dx=\frac12\sum_{k\in \Lambda}|a_k|^2(\mathrm{Id}-k\otimes k).\]

We invoke the following geometric lemma, proven in \cite{BV1, DS1}, which asserts that given a symmetric matrix $R$ close to the identity matrix, we can choose several Beltrami waves $W,$ dependent smoothly on $R,$ such that $$\displaystyle\stackinset{c}{}{c}{}{-\mkern4mu}{\displaystyle\int_{\mathbb{T}^3}} W\otimes W\, \mathrm{d}x = R,$$
which shall be used to construct cancellations in the iteration process.

\begin{Lemma}[Geometric lemma]\label{le-geo}
For any $N\in\mathbb N$, we can find $\varepsilon_\gamma>0$ and $\lambda>1$ with the following property. Let $B_{\varepsilon_\gamma}(\mathrm {Id})$ be the ball in $\mathcal{S}^{3 \times 3}$ centered at $\mathrm{Id}$ with radius $\varepsilon_\gamma$. There exists pairwise disjoint subsets
\[\Lambda_\alpha\subset \mathbb S^2\cap\mathbb Q^3, \ \ \alpha\in\{1, ..., N\},\]
with $\lambda \Lambda_{\alpha}\in \mathbb Z^3$, and smooth positive functions
\[\gamma_{k}^{\alpha}\in C^\infty(B_{\varepsilon_\gamma}(\mathrm{Id})), \ \ \alpha\in\{1, ..., N\}, \ \ k\in\Lambda_\alpha,\]
with derivatives that are bounded independently of $\lambda$, such that -- 
\begin{enumerate}
\item $k\in\Lambda_{\alpha}$ implies $-k\in\Lambda_{\alpha}$ and $\gamma_{k}^\alpha=\gamma_{-k}^\alpha,$
\item For each $R\in B_{\varepsilon_\gamma}(\mathrm{Id})$ we have the identity
\[R=\frac12\sum_{k\in\Lambda_{\alpha}}\left(\gamma_{k}^\alpha(R)\right)^2(\mathrm{Id}-k\otimes k).\]
\end{enumerate}
\end{Lemma}

\subsection{Auxiliary estimates}

We list a few estimates which will be used in upcoming sections. The following $L^p$ de-correlation lemma can be found in \cite{BV, BV1}.
\begin{Lemma}\label{le-cor} 
For given integers $M, \kappa \geq 1,$ let $\lambda\geq 1$ satisfy
\[\frac{2\sqrt 3\pi\lambda}{\kappa}\leq \frac13  \ \ \mbox{and} \ \ \lambda^4\frac{(2\sqrt 3\pi\lambda)^M}{\kappa^M}\leq 1.\]
Let $p\in\{1,2\},$ and let $f$ be a $\mathbb T^3$-periodic function such that 
$$\|D^jf\|_{L^p}\leq C_f \lambda^j, \ \forall j \in [1, M+4]$$
for a certain constant $C_f.$ In addition, let $g$ be a $(\mathbb T/\kappa)^3$-periodic function. Then the following inequality 
\[\|fg\|_{L^p}\lesssim C_f\|g\|_{L^p}\]
holds, where the implicit constant is universal.
\end{Lemma}

We will use the following commutator estimate, proven in \cite{BV1}.
\begin{Lemma}\label{le-comm}
Fix $\kappa\geq 1$, $p\in(1,2]$ and a sufficient large $L \in\mathbb N.$ Let $a\in C^L(\mathbb T^3)$ be such that there exist $\lambda\in [1, \kappa]$ and $C_a>0$ satisfying
$$\|D^ja\|_{L^\infty}\leq C_a\lambda^j, \ \forall j \in [0, L].$$ For $f\in L^p(\mathbb T^3)$ satisfying $\int_{\mathbb T^3}a(x)\mathbb P_{\geq \kappa}f(x)\, \mathrm dx=0$, it holds true that
\[\left\||\nabla|^{-1}(a\mathbb P_{\geq \kappa}f)\right\|_{L^p}\lesssim C_a\left(1+\frac{\lambda^L}{\kappa^{L-2}}\right)\frac{\|f\|_{L^p}}{\kappa},\]
with the implicit constant depending on $p$ and $L$.
\end{Lemma}

In \cite{BV1, DS1}, the following inverse divergence operator was introduced, along with the Calder\'on-Zygmund and Schauder estimates associated with it.
\begin{Lemma}\label{le-anti}
Let $v \in C^\infty(\mathbb T^3)$ be a smooth vector field. There exists a linear operator $\mathcal R$ such that $\mathcal R v(x)$ is a symmetric trace- free matrix for any $x \in \mathbb T^3,$ and
\begin{equation}\notag
\nabla\cdot \mathcal Rv=v-\displaystyle\stackinset{c}{}{c}{}{-\mkern4mu}{\displaystyle\int_{\mathbb T^3}} v\, \mathrm dx.
\end{equation}
In addition, the following estimates on $\mathcal R$ hold for $1<p<\infty.$
\begin{equation}\notag
\|\mathcal R\|_{L^p\to W^{1,p}}\lesssim 1, \ \ \|\mathcal R\|_{C^0\to C^0}\lesssim 1, \ \ 
\|\mathcal R\mathbb P_{\neq0}u\|_{L^p}\lesssim \left\||\nabla|^{-1}\mathbb P_{\neq0}u\right\|_{L^p}.
\end{equation}
\end{Lemma}

\bigskip

\section{Proof of Theorem \ref{thm2}}\label{pt}

To start the iteration process, we set the initial magnetic potential as 
\begin{equation}\notag
A_0(t)=\frac{t}{\lambda_0 (2\pi)^3}\left(0, \cos(\lambda_0x_1), -\sin(\lambda_0x_1)\right).
\end{equation}
We note that here $\left(0, \cos(\lambda_0x_1), -\sin(\lambda_0x_1)\right)=\cos(\lambda_0k_1\cdot x)k_2-\sin(\lambda_0k_1\cdot x)k_3$ with $k_1=(1,0,0)$, $k_2=(0,1,0)$, and $k_3=(0,0,1)$ is a Beltrami wave. We can verify that $B_0=\nabla\times A_0=\lambda_0 A_0$ and $(\nabla\times B_0)\times B_0=0$. Thus $A_0$ satisfies the first iteration
\[\partial_t A_0+(\nabla\times B_0)\times B_0=\nabla\cdot R_0\]
with the symmetric and traceless stress tensor
\begin{equation}\notag
R_0=\frac{1}{\lambda_0^{2}(2\pi)^3}
\begin{bmatrix}
0 & \sin(\lambda_0x_1) & \cos(\lambda_0x_1)\\
\sin(\lambda_0x_1) & 0 & 0\\
\cos(\lambda_0x_1) & 0 & 0\\
\end{bmatrix}.
\end{equation}

A straightforward computation shows that the estimates \eqref{est-induction} holds for $q=0.$ Starting from $(A_0, B_0, R_0)$, we recursively apply Proposition \ref{prop}, resulting in a sequence of approximate solutions $\{(A_q, B_q, R_q)\}_{q \geq 1}$ whose $q$-th element satisfies the corresponding system (\ref{Aq}) and estimates (\ref{est-induction})-(\ref{est-increm}) at the $q$-th level. 

It follows from estimates \eqref{est-increm} that
\begin{equation}\notag
\sum_{q\geq0}\|B_{q+1}-B_q\|_{L^2}= \sum_{q\geq0}\|w_{q+1}\|_{L^2} \lesssim \sum_{q\geq0}\delta_{q+1}^{\frac{1}{2}}\lesssim 1,
\end{equation}
implying the strong convergence of $\{B_q\}_{q=0}^\infty$ to a certain function $B$ in $C^0(0,T;L^2)$. We further infer that $\{A_q\}_{q=0}^\infty$ converges strongly to a function $A$ in $C^0(0,T;H^1)$ satisfying $\nabla\cdot A=0$ and $\nabla\times A=B.$ 

Since $\|R_q\|_{L^\infty(0,T; L^1)}\to 0$ as $q\to \infty$, $(A, B)$ is a weak solution of system \eqref{amhd}. (Equivalently, $B$ is a weak solution of system (\ref{emhd}).) Clearly, $A\in L^\infty(0,T; H^1(\mathbb T^3))$ and $B\in L^\infty(0,T; L^2(\mathbb T^3)).$ 

We denote the magnetic energy of $B$ by
\[\mathcal E(t)=\int_{\T^3}|B(t)|^2\mathrm dx.\]
Similarly, we denote the magnetic energy of $B_0$ by $\mathcal E_0,$ which is given by
\begin{equation}\label{e-e0}
\mathcal E_0(t)=\int_{\T^3}|B_0(t)|^2 \mathrm dx=\frac{ t^2}{(2\pi)^3}.
\end{equation}

The difference $\mathcal E(t)-\mathcal E_0(t)$ is bounded, as
\begin{equation}\notag
\begin{split}
\left|\mathcal E(t)-\mathcal E_0(t)\right|=&\left|\int_{\T^3}\left(|B(t)|^2 - |B_0(t)|^2\right) \mathrm dx\right|\\
\leq & \|B(t)-B_0(t)\|_{L^2}\left(\|B(t)\|_{L^2}+\|B_0(t)\|_{L^2}\right).
\end{split}
\end{equation}
We set $\delta_0$ to be small enough by choosing $a$ and $b$ in \eqref{pp1} large enough, so that $$ \|B-B_0\|_{L^2}\leq\sum_{q\geq0}\|B_{q+1}-B_q\|_{L^2}\leq\sum_{q\geq0} \delta_{q+1}^{\frac12} \leq \frac{1}{(4\pi)^3}.$$
We also note that on the unit time interval $$\|B_0(t)\|_{L^2}\leq (2\pi)^{-3/2} \text{ and } \|B(t)\|_{L^2} \leq 1.$$ Thus, $\|B(t)\|_{L^2}+\|B_0(t)\|_{L^2}\leq 2$ for $t \in [0,1]$. Therefore, we have 
\begin{equation}\label{e-e}
\left|\mathcal E(t)-\mathcal E_0(t)\right|\leq \frac{1}{32 \pi^3}, \ \ t\in[0,1].
\end{equation}

In view of (\ref{e-e0}) and (\ref{e-e}), we obtain
\begin{equation}\notag
\begin{split}
\left|\mathcal E(t)-\mathcal E(0)\right|=&\ \left|\mathcal E(t)-\mathcal E_0(t)+\mathcal E_0(t)-\mathcal E_0(0)+\mathcal E_0(0)-\mathcal E(0)\right|\\
\geq &\ \left|\mathcal E_0(t)-\mathcal E_0(0)\right|-\left|\mathcal E(t)-\mathcal E_0(t)\right|-\left|\mathcal E_0(0)-\mathcal E(0)\right|\\
\geq&\ \frac{t^2}{(2\pi)^3}- \frac{1}{16 \pi^3}.
\end{split}
\end{equation}
It leads to 
\[\left|\mathcal E(1)-\mathcal E(0)\right|\geq \frac{1}{16 \pi^3} >0,\]
which indicates that $\mathcal E(t)$ is not conserved on $[0,1]$.

Analogously, we can check that the magnetic helicity is not conserved. Indeed, at the initial level, we have
\begin{equation}\label{h-h0}
\mathcal H_0(t)=\int_{\T^3}\left( A_0\cdot B_0\right)(t)\mathrm dx=\lambda_0\int_{\T^3}\left(A_0 \cdot A_0\right)(t)\mathrm dx=\frac{t^2}{\lambda_0( 2\pi)^3}.
\end{equation}
The difference between the magnetic helicity $\mathcal H(t)$ of the limit $(A, B)$ and $\mathcal H_0(t)$ satisfy
\[
\begin{split}
\left|\mathcal H(t)-\mathcal H_0(t)\right|=&\ \left|\int_{\T^3}\left(A\cdot B(t)-A_0 \cdot B_0(t)\right)\, dx\right|\\
\leq &\ \|A-A_0\|_{L^2}\|B\|_{L^2}+\|A_0\|_{L^2}\|B-B_0\|_{L^2}.
\end{split}
\]
Note that $\|B\|_{L^2}\leq 1$ thanks to (\ref{est-induction}), $\|A_0\|_{L^2}\leq \lambda_0^{-1}(2\pi)^{-\frac32}$ by construction, and
we can conclude that, for sufficiently small $\delta_0$ (large enough $a$ and $b$)
\[\|A-A_0\|_{L^2}\leq \sum_{q\geq0}\|A_{q+1}-A_q\|_{L^2}\leq \sum_{q\geq0}\lambda_{q+1}^{-1}\delta_{q+1}^{\frac12}\leq \frac{1}{\lambda_0 (4 \pi)^3}\]
\begin{equation}\label{h-h}
\left|\mathcal H(t)-\mathcal H_0(t)\right|\leq \frac{1}{ 32 \pi^3 \lambda_0 }, \ t\in[0,1].
\end{equation}
It follows from (\ref{h-h0}) and (\ref{h-h}) that
\begin{equation}\notag
\begin{split}
\left|\mathcal H(t)-\mathcal H(0)\right|=&\ \left|\mathcal H(t)-\mathcal H_0(t)+\mathcal H_0(t)-\mathcal H_0(0)+\mathcal H_0(0)-\mathcal H(0)\right|\\
\geq &\ \left|\mathcal H_0(t)-\mathcal H_0(0)\right|-\left|\mathcal H(t)-\mathcal H_0(t)\right|-\left|\mathcal H_0(0)-\mathcal H(0)\right|\\
\geq&\ \frac{t^2}{\lambda_0 (2\pi)^3}-\frac{1}{16 \pi^3 \lambda_0 },
\end{split}
\end{equation}
which implies 
\[\left|\mathcal H(1)-\mathcal H(0)\right|\geq \frac{1}{16 \pi^3 \lambda_0}>0.\]
Thus the magnetic helicity $\mathcal H(t)$ is not a constant over $[0,1]$.

\cbdu

\section{Construction of the perturbation} \label{sec-convex}

\subsection{Intermittent Beltrami waves as building blocks} 

We adapt the building blocks which are used to construct non-unique weak solutions for the Navier-Stokes equations in \cite{BV1} as well as for the Hall-MHD system in \cite{D3}, called the intermittent Beltrami waves, to the ideal EMHD system by incorporating only spatial oscillations into the Beltrami waves defined in Section \ref{G}. (In contrast to the construction in \cite{D3}, we exclude temporal oscillations.)

The intermittent Beltrami wave $\mathbb W_{k}$ is defined as
\begin{equation}\label{int-wave}
\mathbb W_{k}(x):=\eta_{k}(x)W_k(x),
\end{equation}
where the oscillation $\eta_{k}(x)$ shall be introduced as follows. For a large integer $r,$ we define the 3D normalized Dirichlet kernel
$$D_r(x):=\frac{1}{(2r+1)^{3/2}}\sum_{k\in\Omega_r}e^{ik\cdot x},$$ where 
$\Omega_r:=\left\{k=(i,j,l): i,j,l\in\{-r,..., r\}\right\}$ is the lattice cube. We fix a small constant $\sigma$ such that $\lambda \sigma \in \mathbb N$ and $\sigma r\ll 1,$ 
and choose an integer $N_0\geq 1$ such that for all $k\in \Lambda_\alpha, \alpha= 1,2,3,..., N,$
\[\{N_0k, N_0k_1, N_0k_2\}\subset \mathbb Z^3.\]
We then define the modified Dirichlet kernel
\begin{equation}\label{eta}
\eta_{k, \lambda,\sigma, r}(x)=
\begin{cases}D_r\left(\lambda\sigma N_0k\cdot x, \lambda\sigma N_0 k_1\cdot x, \lambda\sigma N_0k_2\cdot x\right), \ k\in \Lambda_\alpha^+,\\
\eta_{-k, \lambda,\sigma, r}(x), \ k\in \Lambda_\alpha^-.
\end{cases}
\end{equation}
For simplicity, we use the notation $\eta_{k}(x)=\eta_{k,\lambda,\sigma,r}(x)$. Since $D_r$ satisfies
\begin{equation}\notag
\|D_r\|_{L^2}^2=(2\pi)^3 \text{ and } \|D_r\|_{L^p}\lesssim r^{\frac32-\frac3p}, \ p>1,
\end{equation}
with the implicit constant in the inequality depending only on $p$, we observe that 
\begin{equation}\label{eta-norm}
\displaystyle\stackinset{c}{}{c}{}{-\mkern4mu}{\displaystyle\int_{\mathbb T^3}} \eta_{k}^2(x)\, \mathrm dx = \displaystyle\stackinset{c}{}{c}{}{-\mkern4mu}{\displaystyle\int_{\mathbb T^3}} D_r^2(x)\, \mathrm dx=1 \text{ and }
\|\eta_{k}\|_{L^p}=\|D_{r}\|_{L^p}\lesssim r^{\frac32-\frac3p}, \ p >1.
\end{equation}

Noticing that $\mathbb W_{{k}}$ is supported on certain frequencies, we have 
\begin{gather*}\notag
\mathbb P_{\leq 2\lambda\sigma rN_0}\eta_{{k}}= \eta_{{k}},\\
\mathbb P_{\leq 2\lambda}\mathbb P_{\geq \lambda/2} \mathbb W_{{k}}=\mathbb W_{{k}},\\
\mathbb P_{\leq 4\lambda}\mathbb P_{\geq c_0\lambda} \left( \mathbb W_{{k}}\otimes  \mathbb W_{{k}'}\right)= \mathbb W_{{k}}\otimes  \mathbb W_{{k}'}, \ {k}'\neq -{k}, \ c_0 \text{ a small constant.} 
\end{gather*}
%It is easy to see that
%\begin{equation}\notag
%\begin{split}
%\nabla\cdot \mathbb W_{{k}}=&\ \nabla\eta_{{k}}\cdot W_{{k}},\\
%\nabla\times\mathbb W_{{k}}=&\ \lambda \mathbb W_{{k}}+\nabla\eta_{{k}}\times W_{{k}}.
%\end{split}
%\end{equation}
%Parameters $\lambda, \sigma$, and $r$ will be chosen in an appropriate way such that $\nabla\eta_{{k}}\cdot W_{{k}}$ and $\nabla\eta_{{k}}\times W_{{k}}$ are sufficiently small. 
For $\Lambda_\alpha, \varepsilon_\gamma,$ and $\gamma_{k}$ as in Lemma \ref{le-geo}, we have the following geometric lemma on the intermittent Beltrami waves $\mathbb W_k$, for whose proof we refer readers to \cite{BV1}.
%the following geometric lemma, which is a key ingredient in the construction.
\begin{Lemma}\label{le-geo-2}
Assume that the constants $a_{k}\in\mathbb C$ satisfy $\bar a_{k}=a_{-{k}}$.  The vector field
$$\sum_{\alpha}\sum_{{k}\in\Lambda_\alpha}a_{k}\mathbb W_{{k}}(x)$$
is real valued. Moreover, for each matrix $R\in B_{\varepsilon_\gamma}(\mathrm{Id})$ we have
\begin{equation}\label{geo-id}
\sum_{{k}\in \Lambda_\alpha}\left(\gamma_{k}(R)\right)^2\displaystyle\stackinset{c}{}{c}{}{-\mkern4mu}{\displaystyle\int_{\mathbb T^3}} \mathbb W_{{k}}\otimes \mathbb W_{-{k}}\, \mathrm dx=R.
%\sum_{{k}\in \Lambda_\alpha}\left(\gamma_{k}(R)\right)^2 B_{{k}}\otimes B_{-{k}}=R.
\end{equation}
\end{Lemma}

\subsection{Estimates of building blocks}

%The main purpose of adding the oscillation $\eta_{{k}}$ to the Beltrami waves is to make sure the $L^1$ norm of the waves is significantly smaller than the $L^2$ norm. This can be seen in the following lemma.

\begin{Lemma}\label{le-w} \cite{BV}
The bounds 
\begin{equation}
\|\nabla^m \mathbb W_{{k}}\|_{L^p}\lesssim \lambda^mr^{\frac32-\frac3p}, \label{w-lp}\\
\end{equation}
\begin{equation}
\|\nabla^m\mathbb \eta_{{k}}\|_{L^p}\lesssim (\lambda\sigma r)^mr^{\frac32-\frac3p}  \label{eta-lp}
\end{equation}
hold for all $1<p\leq \infty$.
\end{Lemma}

\begin{Lemma}\label{le-a}
The following bounds hold
\bg\label{a-l2}
\|a_{{k}}\|_{L^2}\lesssim \delta_{q+1}^{\frac12},
\ed
\bg\label{a-lp}
\|a_{{k}}\|_{L^p}\lesssim \delta_{q+1}^{\frac12}\ell^{-\frac12(1-\frac1p)}, \ \ \mbox {for} \ \ p\geq 1,
\ed
\bg\label{a-cn}
\|a_{{k}}\|_{C^m_{x,t}}\lesssim \ell^{-m}, \ \ \mbox {for} \ \ m\geq 1.
\ed
%\bg\label{a-n2}
%\|D^Na_{{k}}\|_{L^2}\lesssim \delta_{q+1}^{\frac12}\ell^{-2N}, \ \ \mbox {for} \ \ N\geq 1,
%\ed
%\bg\label{a-n-infty}
%\|D^Na_{{k}}\|_{L^\infty}\lesssim \delta_{q+1}^{\frac12}\ell^{-2N-\frac12}, \ \ \mbox {for} \ \ N\geq 1.
%\ed
\end{Lemma}

\subsection{Mollification and amplitude functions}

To avoid the loss of derivatives, we need to mollify the constructed solutions. The mollification can be done through a standard procedure by using Friedrichs mollifiers. We define the mollified versions of $A_q$, $B_q$, and $R_q$ in space and time at the scale $\ell$ as
\[A_\ell=(A_q*_x\phi_{\ell})*_t\varphi_\ell, \ B_\ell=(B_q*_x\phi_{\ell})*_t\varphi_\ell, \ R_\ell=(R_q*_x\phi_{\ell})*_t\varphi_\ell,\]
where $\phi_\ell$ and $\varphi_\ell$ are standard Friedrichs mollifiers on $\mathbb R^3$ and $\mathbb R,$ respectively. The mollified triple $(A_\ell, B_\ell, R_\ell)$ satisfies
\begin{equation}\notag
\begin{split}
\partial_t A_\ell+\nabla\cdot (B_\ell\otimes B_\ell)- \nabla p_\ell =\nabla\cdot (R_\ell+R_{\mathrm{comm}}),\\
\nabla\cdot A_\ell=\nabla\cdot B_\ell=0,
\end{split}
\end{equation}
with the traceless symmetric commutator stress tensor $R_{\mathrm{comm}}$ defined as
\[R_{\mathrm{comm}}=(B_{\ell} \mathring\otimes B_{\ell})-\left((B_q \mathring\otimes B_q)*_x\phi_\ell\right)*_t\varphi_\ell,\]
and the pressure term $p_\ell$ defined as
\[p_\ell=(p_q*_x\phi_\ell)*_t\varphi_\ell-|B_\ell|^2+(|B_q|^2*_x\phi_\ell)*_t\varphi_\ell.\]
We have the following estimates on the stress tensors --
\begin{equation}\label{Rl1}
\|\nabla^mR_\ell\|_{L^1}\lesssim \ell^{-m}\|R_q\|_{L^1}, \ \|R_{\mathrm{comm}}\|_{L^1}\lesssim \|R_{\mathrm{comm}}\|_{C^0} \lesssim \ell^2\|B_q\|_{C^1_{x,t}}^2.
\end{equation}

%\subsection{Amplitude functions}

We prescribe appropriate amplitude functions in view of the geometric lemma. To this end, we choose a smooth function
$\chi:[0,\infty)\to \mathbb R,$ defined as 
\begin{equation}\notag
\chi(z)=
\begin{cases}
1,  \ 0\leq z\leq 1,\\
z, \ z\geq 2,
\end{cases}
\end{equation}
with $1 \leq  \chi(z) \leq 2$ and $0 \leq \chi'(z)\leq 1$ for $z \in [1,2]$. We define
\[\rho(x,t)=2\delta_{q+1}\varepsilon^{-1}c\chi\left((c\delta_{q+1})^{-1}|R_\ell (x,t)|\right).\]
Notice that we have the following properties of $\rho$:
\[\left|\frac{R_{\ell}(x,t)}{\rho(x,t)}\right|\leq \varepsilon,\]
\[\|\rho\|_{L^p}\leq 8\varepsilon^{-1}\left(c(8\pi^3)^{\frac1p}\delta_{q+1}+\|R_\ell\|_{L^p}\right), \ \ \ p\in[1,\infty),\]
\[\|\rho\|_{C^0_{x,t}}\lesssim \ell^{-3}, \ \ \ \|\rho\|_{C^j_{x,t}}\lesssim \ell^{-4j}, \ \ \ j\geq 1,\]
\[\|\rho^{\frac{1}{2}}\|_{C^0_{x,t}}\lesssim \ell^{-2}, \ \ \ \|\rho^{\frac{1}{2}}\|_{C^j_{x,t}}\lesssim \ell^{-5j}, \ \ \ j\geq 1.\]
The amplitude functions $a_{(k)}$ are defined as
\begin{equation}\label{def-ak}
a_{(k)}(x,t)=\rho^{\frac{1}{2}}\gamma_{(k)}\left(\mathrm{Id}-\frac{R_{\ell}}{\rho}\right), \ \ \ k\in\Lambda.
\end{equation}
We have the estimates for $a_{(k)}$,
\begin{equation}\notag
\begin{split}
\|a_{(k)}\|_{L^2}\lesssim \delta_{q+1}^{\frac12}, \  \|a
_{(k)}\|_{C^j_{x,t}}\lesssim \ell^{-5j-2}.
\end{split}
\end{equation}
In view of (\ref{def-ak}) and Lemma \ref{le-geo-2}, we conclude
\begin{equation}\label{Rq-cancel}
\sum_{k+k'=0}a_k^2 \displaystyle\stackinset{c}{}{c}{}{-\mkern4mu}{\displaystyle\int_{\mathbb T^3}} \mathbb W_{{k}}\otimes \mathbb W_{k'}\, \mathrm dx=\rho \mathrm{Id}-R_\ell.
%+\sum_{k\in \Lambda}a_k^2\mathbb P_{\neq 0}(\eta_k^2)W_k\otimes W_k.
\end{equation}

We introduce the the following parameters for $q \in \mathbb{N}.$

\begin{equation}\label{pp2}
r=\lambda_{q+1}^{\frac23}, \ \ \sigma=\lambda_{q+1}^{-\frac{5}{6}},  \ \ \ell=\lambda_q^{-20},
\end{equation}
where the constants $a\gg 1$ and $b\gg 1$ are chosen to be sufficiently large and the constant $\beta\ll 1$ is positive and sufficiently small.

\subsection{Construction of the perturbation}
\label{sec-design}
The perturbation $v_{q+1}=A_{q+1}-A_q$ consists of a principal part and a corrector --
$$ v_{q+1} :=v_{q+1}^{p}+v_{q+1}^{c},$$
where $v_{q+1}^{p}$ and $v_{q+1}^{c}$  are defined as
\begin{equation}\notag
\begin{split}
v_{q+1}^{p}:=&\sum_{{k}\in \Lambda_\alpha} a_{k}\mathbb W_{{k}}=\lambda_{q+1}^{-1}\sum_{{k}\in \Lambda_\alpha}a_{{k}}\eta_{{k}} W_{{k}},\\
v_{q+1}^{c}:=& \lambda_{q+1}^{-2}\sum_{{k}\in \Lambda_\alpha}\nabla(a_{{k}}\eta_{{k}})\times W_{{k}}.\\
%v_{q+1}^{t}=&\mu^{-1}\lambda_{q+1}^{-1}\sum_{{k}}\mathbb P_H\mathbb P_{\neq0}(a_{k}^2\eta^2{k}).
\end{split}
\end{equation}
We can verify that
\[\nabla\cdot (v_{q+1}^p+v_{q+1}^c)=\lambda_{q+1}^{-2}\sum_{{k}\in \Lambda_\alpha}\nabla\cdot(\nabla\times (a_{{k}}\eta_{{k}}W_{{k}}))=0.\]

In view of the construction of $v_{q+1},$ we define the perturbation $w_{q+1}=B_{q+1}-B_q$ as
$$w_{q+1}:=w_{q+1}^p+w_{q+1}^c:= \nabla\times v_{q+1}^{p}+\nabla\times v_{q+1}^{c}.$$
It is clear that $\nabla\cdot w_{q+1}^p=\nabla\cdot w_{q+1}^c=\nabla\cdot w_{q+1}=0$ and $w_{q+1}=\nabla\times v_{q+1}.$

\section{Proof of Proposition \ref{prop}}
\label{sec-prop}
To close the analysis of the previous section, we only need to prove Proposition \ref{prop}, which will be achieved through a series of estimates.

\subsection{Estimates of the perturbation}

\begin{Lemma}\label{le-vq}
The increment $v_{q+1}=A_{q+1}-A_q$ satisfies the following estimates
\begin{equation}\label{v-p2}
\|v_{q+1}^p\|_{L^2}\lesssim \lambda_{q+1}^{-1}\delta_{q+1}^{\frac12},
\end{equation}
\begin{equation}\label{v-c2}
\|v_{q+1}^c\|_{L^2}\lesssim \ell^{-1}\lambda_{q+1}^{-1}\delta_{q+1}^{\frac12},
\end{equation}
\begin{equation}\label{v-pp}
\|v_{q+1}^p\|_{L^p}\lesssim \lambda_{q+1}^{-1}\delta_{q+1}^{\frac12}\ell^{-\frac12(1-\frac1p)}r^{\frac32-\frac3p}, \ \ p\geq1,
\end{equation}
\bg\label{v-cp}
\|v_{q+1}^c\|_{L^p}\lesssim \lambda_{q+1}^{-1}\delta_{q+1}^{\frac12}\ell^{-\frac12(1-\frac1p)}\sigma r^{\frac52-\frac3p}, \ \ p\geq 1,
\ed
%\begin{equation}\label{v-w1p}
%\|v_{q+1}^p\|_{W^{1,p}}+\|v_{q+1}^c\|_{W^{1,p}}\lesssim \ell^{-2}r^{\frac32-\frac3p}, \ \ p\geq 1,
%\end{equation}
\begin{equation}\label{v-cn}
\|v_{q+1}^p\|_{C^m_{x,t}}+\|v_{q+1}^c\|_{C^m_{x,t}}\lesssim \lambda_{q+1}^{\frac{1+5m}{2}},
\end{equation}
\begin{equation}\label{b-cn}
\|A_{q+1}\|_{C^m_{x,t}}\lesssim \lambda_{q+1}^{\frac{3+5m}{2}}.
%\ \ \||\nabla|^NB_{q+1}\|_{L^p}\lesssim \lambda_{q+1}^{N-1}r^{\frac32-\frac3p}.
\end{equation}
\end{Lemma}

\begin{Lemma}\label{le-w}
The increment $w_{q+1}=B_{q+1}-B_q$ satisfies the following estimates,
\bg\label{wp-l2}
\|w_{q+1}^p\|_{L^2}\lesssim \delta_{q+1}^{\frac12},
\ed
\begin{equation}\label{wc-l2}
\|w_{q+1}^c\|_{L^2}\lesssim \ell^{-1}\delta_{q+1}^{\frac12},
\end{equation}
\bg\label{w-p}
\|w_{q+1}\|_{L^p}\lesssim \delta_{q+1}^{\frac12}\ell^{-\frac12(1-\frac1p)}r^{\frac32-\frac3p}, \ \ p\geq 1,
\ed
\bg\label{w-cp}
\|w_{q+1}^c\|_{L^p}\lesssim \delta_{q+1}^{\frac12}\ell^{-\frac12(1-\frac1p)}\sigma r^{\frac52-\frac3p}, \ \ p\geq 1,
\ed
\begin{equation}\label{w-1p}
\|w_{q+1}^p\|_{W^{1,p}}+\|w_{q+1}^c\|_{W^{1,p}}\lesssim \ell^{-2}\lambda_{q+1}r^{\frac32-\frac3p}, \ \ p\geq 1,
\end{equation}
%\begin{equation}\label{wt-p}
%\|\partial_tw_{q+1}^p\|_{L^p}+\|\partial_tw_{q+1}^c\|_{L^p}\lesssim \ell^{-2}\lambda_{q+1}\sigma\mu r^{\frac52-\frac3p}, \ \ p\geq 1,
%\end{equation}
%\begin{equation}\label{w-cn}
%\|w_{q+1}^p\|_{C^m_{x,t}}+\|w_{q+1}^c\|_{C^m_{x,t}}\lesssim \lambda_{q+1}^{\frac{3+5m}{2}},
%\end{equation}
%\begin{equation}\label{wn-p}
%\||\nabla|^mw_{q+1}^p\|_{L^p}+\||\nabla|^mw_{q+1}^c\|_{L^p}\lesssim \lambda_{q+1}^{m} r^{\frac32-\frac3p}, \ \ p\geq 1,
%\end{equation}
\bg\label{Jn-p}
\||\nabla|^m B_{q+1}\|_{L^p}\lesssim \lambda_{q+1}^{m+\frac32},  \ \ p\geq 1.
\ed
\end{Lemma}

The conclusion (\ref{est-increm}) follows directly from (\ref{wp-l2}), (\ref{wc-l2}) and (\ref{w-1p}). We are only left to estimate the stress tensor at level $q+1$, which is the main task of next subsection.

\subsection{Estimate of the stress tensor $R_{q+1}$}
\label{sec-R}

\begin{Lemma}\label{le-R}
Consider the equation
\begin{equation}\label{eq-RM}
\partial_tA_{q+1}+\nabla\cdot (B_{q+1}\otimes B_{q+1})+\nabla p_{q+1}=\nabla\cdot  R_{q+1}.
\end{equation}
There exists a traceless symmetric tensor $R_{q+1}$ satisfying (\ref{eq-RM}). 
In addition, there exists $p>1$ sufficiently close to 1, and a sufficiently small $\varepsilon_R>0$ independent of $q$ such that 
\begin{equation}\label{est-RM}
\|R_{q+1}\|_{L^p}\lesssim \lambda_{q+1}^{-2\varepsilon_R}\delta_{q+2}
%\|\tilde M_{q+1}\|_{L^p}\lesssim \lambda_{q+1}^{-2\varepsilon_R}\delta_{q+2}
\end{equation}
holds for some implicit constant which depends on $p$ and $\varepsilon_R$.
\end{Lemma}
\pf
%To estimate the stress tensor $R_{q+1}$, 
We first subtract the equation (\ref{Aq}) at level of $A_q$ from the equation at level of $A_{q+1}$ to arrive
%\begin{equation}\notag
%\partial_t v_{q+1}+\nabla\cdot (B_{q+1}\otimes B_{q+1}-B_q\otimes B_q)-\nabla \tilde p_{q+1}
%=\nabla\cdot R_{q+1}-\nabla\cdot R_q.
%\end{equation}
%Rearranging the terms we obtain
\begin{equation}\notag
\begin{split}
\nabla\cdot R_{q+1}=&\ \partial _t v_{q+1}+\nabla\cdot (B_{q+1}\otimes w_{q+1}+ w_{q+1}\otimes B_q)+\nabla\cdot R_q-\nabla \tilde p_{q+1}\\
=&\ \nabla\cdot[\mathcal R\left(\partial _t v_{q+1}^p+\partial _t v_{q+1}^c\right)+ (B_q\otimes w_{q+1}+w_{q+1}\otimes B_q)]\\
&+\nabla\cdot (w_{q+1}^c\otimes w_{q+1}+w_{q+1}^p\otimes w_{q+1}^c)\\
&+\nabla\cdot (w_{q+1}^p\otimes w_{q+1}^p+R_q) +\nabla \tilde p_{q+1}\\
=&: \nabla\cdot R_{\mathrm{linear}}+\nabla\cdot R_{\mathrm{Nash}}+\nabla\cdot R_{\mathrm{oscillation}}+\nabla \tilde p_{q+1}.
\end{split}
\end{equation}
%We further classify the terms on the right hand side into linear, correction and oscillation terms:
%On the right hand side of the equation above, the first three lines correspond to linear terms, the middle three correspond to correction terms, and the last two lines correspond to oscillation terms. The estimates of them will be accomplished separately below.
We estimate the linear error $R_{\mathrm{linear}}$, Nash error $R_{\mathrm{Nash}}$, and oscillation error $R_{\mathrm{oscillation}}$ separately in the following. 

Recall that $v_{q+1}^p=\lambda_{q+1}^{-1}\sum_{k\in \Lambda_\alpha}a_k\eta_k W_k$ and $v_{q+1}^c=\lambda_{q+1}^{-2}\sum_{k\in \Lambda_\alpha}\nabla(a_k\eta_k)\times W_k$, we have
\[v_{q+1}^p+v_{q+1}^c=\lambda_{q+1}^{-2}\sum_{k\in \Lambda_\alpha}\nabla\times(a_k\eta_k W_k),\]
and hence
\[\partial_tv_{q+1}^p+ \partial_tv_{q+1}^c=\lambda_{q+1}^{-2}\sum_{k\in \Lambda_\alpha}\nabla\times (\partial_ta_k\eta_kW_k).\]
It follows 
\begin{equation}\notag
\begin{split}
\|\mathcal R\left(\partial_tv_{q+1}^p+ \partial_tv_{q+1}^c\right)\|_{L^p}\lesssim& \lambda_{q+1}^{-2}\sum_{k\in \Lambda_\alpha}\|\partial_ta_k\eta_kW_k\|_{L^p}\\
\lesssim &\lambda_{q+1}^{-2}\sum_{k\in \Lambda_\alpha}\|a_k\|_{C^1_{x,t}}\|\eta_kW_k\|_{L^p}\\
\lesssim &\lambda_{q+1}^{-2}\ell^{-1}r^{\frac32-\frac3p}
\end{split}
\end{equation}
where we used (\ref{w-lp}) and (\ref{a-cn}) in the last step. By (\ref{b-cn}) and (\ref{w-p}), we have
\begin{equation}\notag
\begin{split}
&\|B_{q}\otimes w_{q+1}+ w_{q+1}\otimes B_{q}\|_{L^p}\\
\lesssim &\|B_q\|_{L^\infty}\|w_{q+1}\|_{L^p}\\
\lesssim &\ \lambda_q^4\ell^{-\frac12(1-\frac1p)}\delta_{q+1}^{\frac12}r^{\frac32-\frac3p}.
\end{split}
\end{equation} 
The last two estimates imply that 
\begin{equation}\label{est-RL}
\|R_{\mathrm{linear}}\|_{L^p}\lesssim \lambda_{q+1}^{-2}\ell^{-1}r^{\frac32-\frac3p}+\lambda_q^4\ell^{-\frac12(1-\frac1p)}\delta_{q+1}^{\frac12}r^{\frac32-\frac3p}.
\end{equation}

%Estimate of $\|w_{q+1}^c\|_{L^p}$: 
%\[w_{q+1}^c=\nabla\times v_{q+1}^c=\lambda_{q+1}^{-2}\sum \nabla\times(\nabla(a_k\eta_k)\times W_k)=\lambda_{q+1}^{-2}\sum[-W_k\nabla\cdot\nabla(a_k\eta_k)+W_k\cdot \nabla\nabla(a_k\eta_k)-\nabla(a_k\eta_k)\cdot\nabla W_k]\]
%\[\|w_{q+1}^c\|_{L^p}\lesssim \lambda_{q+1}^{-2}\left(\sum \|\nabla\nabla(a_k\eta_k)\|_{L^p}+\lambda_{q+1}\sum \|\nabla(a_k\eta_k)\|_{L^p}\right)\]
%\[\|\nabla\nabla(a_k\eta_k)\|_{L^p}\lesssim \|a_k\|_{C^2_{x,t}}(\|\eta_k\|_{L^p}+\|\nabla\eta_k\|_{L^p}+\|\nabla^2\eta_k\|_{L^p})\lesssim \ell^{-2}(1+\lambda\sigma r+(\lambda\sigma r)^2)r^{\frac32-\frac3p}\]
%\[\|\nabla(a_k\eta_k)\|_{L^p}\lesssim \|a_k\|_{C^_{x,t}}(\|\eta_k\|_{L^p}+\|\nabla\eta_k\|_{L^p})\lesssim \ell^{-1}(1+\lambda\sigma r)r^{\frac32-\frac3p}\]
% \[\|w_{q+1}^c\|_{L^p}\lesssim \lambda_{q+1}^{-2}\ell^{-2}r^{\frac32-\frac3p}(1+\lambda\sigma r+(\lambda\sigma r)^2)+\lambda_{q+1}^{-1}\ell^{-1}r^{\frac32-\frac3p}(1+\lambda\sigma r)\]

The Nash error can be estimated as, by using (\ref{w-p}) and (\ref{w-cp}), 
\begin{equation}\label{est-RN}
\begin{split}
\|R_{\mathrm{Nash}}\|_{L^p}\leq&\ 2\|w_{q+1}\|_{L^{2p}}\|w_{q+1}^c\|_{L^{2p}}\\
\lesssim&\ \delta_{q+1}\ell^{-(1-\frac1p)}\sigma r^{4-\frac3p}.\\
%\lesssim &\ \lambda_{q+1}^{-2\varepsilon_R}\delta_{q+1}.
\end{split}
\end{equation}

Next we estimate the oscillation error. Recall 
\begin{equation}\notag
\begin{split}
w_{q+1}^p=&\ \nabla\times v_{q+1}^p\\
=&\ \lambda_{q+1}^{-1}\sum_{k\in \Lambda_\alpha}\nabla(a_k\eta_k)\times W_k+\sum_{k\in \Lambda_\alpha}a_k\eta_kW_k\\
=&: W_\epsilon+\sum_{k\in \Lambda_\alpha}a_k\eta_kW_k.
\end{split}
\end{equation}
It follows 
\begin{equation}\notag
\begin{split}
\nabla\cdot R_{\mathrm{oscillation}}=&\ \nabla\cdot(w_{q+1}^p\otimes w_{q+1}^p)+\nabla\cdot R_q\\
=&\ \nabla\cdot\left(\sum_{k+k'=0}a_ka_{k'}\eta_k\eta_{k'}W_k\otimes W_{k'}+R_q\right)\\
&+\nabla\cdot\left(\sum_{k+k'\neq 0}a_ka_{k'}\eta_k\eta_{k'}W_k\otimes W_{k'}\right)+\nabla\cdot (W_\epsilon \otimes W_\epsilon)\\
&+\nabla\cdot\left(W_\epsilon\otimes \left(\sum_{k\in \Lambda_\alpha}a_k\eta_kW_k\right)\right)+\nabla\cdot\left(\left(\sum_{k\in \Lambda_\alpha}a_k\eta_kW_k\right) \otimes W_\epsilon\right)\\
=&: \nabla\cdot R_{\mathrm{o,1}}+\nabla\cdot R_{\mathrm{o,2}}+\nabla\cdot R_{\mathrm{o,3}}+\nabla\cdot R_{\mathrm{o,4}}+\nabla\cdot R_{\mathrm{o,5}}.
\end{split}
\end{equation}
Appealing (\ref{Rq-cancel}), we obtain 
\begin{equation}\notag
\begin{split}
\nabla\cdot R_{\mathrm{o,1}}=&\ \sum_{k+k'=0}\nabla\cdot\left(a_ka_{k'}\left(\mathbb W_k\otimes\mathbb W_{k'}-\displaystyle\stackinset{c}{}{c}{}{-\mkern4mu}{\displaystyle\int_{\mathbb T^3}} \mathbb W_{{k}}\otimes \mathbb W_{k'}\, \mathrm dx\right)\right)+\nabla \rho\\
=&\ \sum_{k+k'=0}\nabla\cdot\left(a_ka_{k'}\mathbb P_{\geq \frac12\lambda_{q+1}\sigma}(\mathbb W_k\otimes \mathbb W_{k'})\right)+\nabla \rho\\
=&\ \sum_{k\in\Lambda}\mathbb P_{\neq 0}\left(\mathbb P_{\geq \frac12\lambda_{q+1}\sigma}(\mathbb W_k\otimes \mathbb W_{-k}+\mathbb W_{-k}\otimes \mathbb W_{k})\nabla a_k^2\right)\\
&+\sum_{k\in\Lambda}\mathbb P_{\neq 0}\left(a_k^2\nabla\cdot(\mathbb W_k\otimes \mathbb W_{-k}+\mathbb W_{-k}\otimes \mathbb W_{k})\right)+\nabla \rho\\
=&\ \sum_{k\in\Lambda}\mathbb P_{\neq 0}\left(\mathbb P_{\geq \frac12\lambda_{q+1}\sigma}(\mathbb W_k\otimes \mathbb W_{-k}+\mathbb W_{-k}\otimes \mathbb W_{k})\nabla a_k^2\right)\\
&+\sum_{k\in\Lambda}\mathbb P_{\neq 0}\left(a_k^2\nabla \eta_k^2\right)+\nabla \rho\\
=&\ \sum_{k\in\Lambda}\mathbb P_{\neq 0}\left(\mathbb P_{\geq \frac12\lambda_{q+1}\sigma}(\mathbb W_k\otimes \mathbb W_{-k}+\mathbb W_{-k}\otimes \mathbb W_{k})\nabla a_k^2\right)\\
&-\sum_{k\in\Lambda}\mathbb P_{\neq 0}\left(\mathbb P_{\geq\frac12\lambda_{q+1}\sigma}(\eta_k^2)\nabla a_k^2\right)+\sum_{k\in\Lambda}\nabla\left(a_k^2\mathbb P_{\geq\frac12\lambda_{q+1}\sigma}(\eta_k^2)\right)+\nabla \rho.\\
\end{split}
\end{equation}
Neglecting the pressure terms in the equation above, $R_{\mathrm{o,1}}$ can be estimated by using (\ref{a-cn}) and the commutator estimate (\ref{le-comm}), for large enough $m$
\begin{equation}\notag%\label{est-RO1}
\begin{split}
\|R_{\mathrm{o,1}}\|_{L^p}\lesssim&\ \sum_{k\in \Lambda}\|\mathcal R\mathbb P_{\neq 0}\left(\mathbb P_{\geq \frac12\lambda_{q+1}\sigma}(\mathbb W_k\otimes \mathbb W_{-k}+\mathbb W_{-k}\otimes \mathbb W_{k})\nabla a_k^2\right) \|_{L^p}\\
&+\sum_{k\in \Lambda}\|\mathcal R \mathbb P_{\neq 0}\left(\mathbb P_{\geq\frac12\lambda_{q+1}\sigma}(\eta_k^2)\nabla a_k^2\right)\|_{L^p}\\
\lesssim&\ \frac{1}{\ell^2\lambda_{q+1}\sigma}\left(1+\frac{1}{\ell^m(\lambda_{q+1}\sigma)^{m-2}}\right)\sum_{k\in\Lambda}\left(\|\mathbb W_k\otimes \mathbb W_{-k}\|_{L^p}+\|\eta_k^2\|_{L^p}\right)\\
\lesssim&\ \frac{1}{\ell^2\lambda_{q+1}\sigma}\sum_{k\in\Lambda}\left(\|\mathbb W_k\|_{L^{2p}}^2+\|\eta_k\|_{L^{2p}}^2\right)\\
\lesssim&\ \frac{r^{3-\frac3p}}{\ell^2\lambda_{q+1}\sigma}.
\end{split}
\end{equation}
Regarding $R_{\mathrm{o,2}}$, we have
\begin{equation}\notag
\begin{split}
\nabla\cdot R_{\mathrm{o,2}}
=&\ \sum_{k+k'\neq0}\nabla\cdot\left(a_ka_{k'}\mathbb P_{\geq c\lambda_{q+1}}(\mathbb W_k\otimes \mathbb W_{k'})\right)\\
=&\ \sum_{k+k'\neq0}\mathbb P_{\neq 0}\left(\mathbb P_{\geq c\lambda_{q+1}}(\mathbb W_k\otimes \mathbb W_{k'})\nabla (a_ka_{k'})\right)\\
&+\frac12\sum_{k+k'\neq0}\mathbb P_{\neq 0}\left(a_ka_{k'}\nabla\cdot(\mathbb W_k\otimes \mathbb W_{k'}+\mathbb W_{k'}\otimes \mathbb W_{k})\right)\\
=&\ \sum_{k+k'\neq0}\mathbb P_{\neq 0}\left(\mathbb P_{\geq c\lambda_{q+1}}(\mathbb W_k\otimes \mathbb W_{k'})\nabla (a_ka_{k'})\right)\\
&+\frac12\sum_{k+k'\neq0}a_ka_{k'}\mathbb P_{\geq c\lambda_{q+1}}\left(\nabla(\eta_k\eta_{k'})(W_k\otimes W_{k'}+W_{k'}\otimes W_k-(W_k\cdot W_{k'})\mathrm{Id}\right)\\
&-\frac12\sum_{k+k'\neq0}\nabla(a_ka_{k'})P_{\geq c\lambda_{q+1}}\left(\mathbb W_k\cdot\mathbb W_{k'}\right)+\frac12\sum_{k+k'\neq0}\nabla(a_ka_{k'}\mathbb W_k\cdot\mathbb W_{k'}).
\end{split}
\end{equation}
Similarly, neglecting the last pressure term, $R_{\mathrm{o,2}}$ can be estimated as
\begin{equation}\notag
\begin{split}
&\|R_{\mathrm{o,2}}\|_{L^p}\\
\lesssim& \sum_{k+k'\neq 0}\|\mathcal R\mathbb P_{\neq 0}\left(\mathbb P_{\geq c\lambda_{q+1}}(\mathbb W_k\otimes \mathbb W_{k'})\nabla (a_ka_{k'})\right)\|_{L^p}\\
&+\sum_{k+k'\neq 0}\|\mathcal Ra_ka_{k'}\mathbb P_{\geq c\lambda_{q+1}}\left(\nabla(\eta_k\eta_{k'})(W_k\otimes W_{k'}+W_{k'}\otimes W_k-(W_k\cdot W_{k'})\mathrm{Id}\right)\|_{L^p}\\
&+\sum_{k+k'\neq 0}\|\mathcal R\nabla(a_ka_{k'})P_{\geq c\lambda_{q+1}}\left(\mathbb W_k\cdot\mathbb W_{k'}\right)\|_{L^p}\\
\lesssim & \frac{1}{\ell^2\lambda_{q+1}}\left(1+\frac{1}{\ell^m\lambda_{q+1}^{m-2}}\right)\sum_{k+k'\neq 0}\left(\|\mathbb W_k\otimes \mathbb W_{k'}\|_{L^p}+\|\nabla(\eta_k\eta_{k'})\|_{L^p}+\|\mathbb W_k\cdot \mathbb W_{k'}\|_{L^p}\right)\\
\lesssim & \frac{1}{\ell^2\lambda_{q+1}}\left(1+\frac{1}{\ell^m\lambda_{q+1}^{m-2}}\right)\left(r^{3-\frac3p}+\lambda_{q+1}\sigma r^{4-\frac3p}\right)\\
\lesssim & \frac{1}{\ell^2\lambda_{q+1}}(1+\lambda_{q+1}\sigma r)r^{3-\frac3p}.\\
\end{split}
\end{equation}
for large enough $m>0$.

The estimates of the rest terms in the oscillation error are trivial. Applying (\ref{a-lp}), (\ref{a-cn}), and (\ref{eta-lp}) leads to
\begin{equation}\notag
\begin{split}
\|R_{\mathrm{o,3}}\|_{L^p}\lesssim&\ \|W_\epsilon\|_{L^{2p}}^2\\
\lesssim&\ \lambda_{q+1}^{-2}\left(\sum_{k\in\Lambda}\|\nabla(a_k\eta_k)\|_{L^{2p}}\right)^2\\
\lesssim&\ \lambda_{q+1}^{-2}\left(\sum_{k\in\Lambda}\|a_k\|_{C^1_{x,t}}\|\eta_k\|_{L^{2p}}+\|a_k\|_{L^\infty}\|\nabla\eta_k\|_{L^{2p}}\right)^2\\
\lesssim&\ \lambda_{q+1}^{-2}\left(\ell^{-1}r^{\frac32-\frac3{2p}}+\delta_{q+1}^{\frac12}\ell^{-\frac12}\lambda_{q+1}\sigma rr^{\frac32-\frac3{2p}}\right)^2\\
\lesssim&\ \lambda_{q+1}^{-2}\ell^{-2}r^{3-\frac3{p}}+\lambda_{q+1}^{-2}\delta_{q+1}\ell^{-1}(\lambda_{q+1}\sigma r)^2r^{3-\frac3{p}}.
\end{split}
\end{equation}
Similarly, we have
\begin{equation}\notag
\begin{split}
&\|R_{\mathrm{o,4}}\|_{L^p}+\|R_{\mathrm{o,5}}\|_{L^p}\\
\lesssim &\ \|W_\epsilon\otimes \sum_{k\in\Lambda}a_k\eta_k W_k\|_{L^p}\\
\lesssim &\ \|W_\epsilon\|_{L^{2p}} \|\sum_{k\in\Lambda}a_k\eta_k W_k\|_{L^{2p}}\\
\lesssim &\ \|W_\epsilon\|_{L^{2p}} \sum_{k\in\Lambda}\|a_k\|_{L^\infty}\|\eta_k\|_{L^{2p}}\\
\lesssim &\ \lambda_{q+1}^{-1}\left(\ell^{-1}r^{\frac32-\frac3{2p}}+\delta_{q+1}^{\frac12}\ell^{-\frac12}\lambda_{q+1}\sigma rr^{\frac32-\frac3{2p}}\right)\delta_{q+1}^{\frac12}\ell^{-\frac12}r^{\frac32-\frac3{2p}}\\
\lesssim &\ \lambda_{q+1}^{-1}\delta_{q+1}^{\frac12}\ell^{-\frac12}\left(\ell^{-1}r^{3-\frac3{p}}+\delta_{q+1}^{\frac12}\ell^{-\frac12}\lambda_{q+1}\sigma rr^{3-\frac3{p}}\right).
\end{split}
\end{equation}
Combining the estimates for $R_{\mathrm{o,i}}$ with $1\leq i\leq 5$, the oscillation error (module the pressure terms) satisfies
\begin{equation}\label{est-RO}
\begin{split}
\|R_{\mathrm{oscillation}}\|_{L^p}\lesssim &\ \frac{r^{3-\frac3p}}{\ell^2\lambda_{q+1}\sigma}+\frac{1}{\ell^2\lambda_{q+1}}(1+\lambda_{q+1}\sigma r)r^{3-\frac3p}\\
&+\lambda_{q+1}^{-2}\ell^{-2}r^{3-\frac3{p}}+\lambda_{q+1}^{-2}\delta_{q+1}\ell^{-1}(\lambda_{q+1}\sigma r)^2r^{3-\frac3{p}}\\
&+\lambda_{q+1}^{-1}\delta_{q+1}^{\frac12}\ell^{-\frac12}\left(\ell^{-1}r^{3-\frac3{p}}+\delta_{q+1}^{\frac12}\ell^{-\frac12}\lambda_{q+1}\sigma rr^{3-\frac3{p}}\right).
\end{split}
\end{equation}
The choice of the parameters (\ref{pp1})-(\ref{pp2}) combined with (\ref{est-RL})-(\ref{est-RO}) leads to the conclusion of the lemma.

\cbdu

%\begin{Remark}In fact, for 3D ideal MHD, the construction of B-Buckmaster-Vicol produced the finite energy weak solutions for which $\mathcal M(t)$ is not conserved.
%\end{Remark}

%{\color{blue}Question: Verify the construction in my paper for the 3D Hall-MHD (non-uniqueness paper), whether $M(t)$ is conserved or not for the weak solutions in Leray-Hopf class or in higher regularity space if $\mu=0$? As the first step, try to verify it for the 3D EMHD $B_t+\nabla\times((\nabla\times B)\times B)=\mu\Delta B$. For $\mu=0$, we can just rescale the construction of the perturbation $w_{q+1}=A_{q+1}-A_q$ by a multiple of $\lambda_{q+1}$.}

%\Endrefs
\end{document}